\documentclass[9pt,amstex]{article}
\usepackage{amssymb,amsmath} 
\usepackage{latexsym}

\newtheorem{dfn}{Definition}[section] 
\newtheorem{rmk}{Remark}[section]
\newtheorem{thm}{Theorem}[section] 
\newtheorem{cor}{Corollary}[section]
\newtheorem{prop}{Proposition}[section] 
\newtheorem{lem}{Lemma}[section]

\newtheorem{ex}{Example}[section] 

\def\cyclic{\mathop{\kern0.9ex{{+}
\kern-2.2ex\raise-.28ex\hbox{\Large\hbox
{$\circlearrowright$}}}}\limits}

\def\buildrel#1_#2^#3{\mathrel{\mathop{\kern 0pt#1}\limits_{#2}^{#3}}}

\newcommand{\Pf}{{\em Proof}. }
\newcommand{\EPf}
{%
\mbox{}%
\nolinebreak%
\hfill%
\rule{2mm}{2mm}%
\medbreak%
\par%
}

\newcommand{\id}{\mbox{$\mathtt{Id}$}}

\newcommand{\Ad}{\mbox{$\mathtt{Ad}$}}

\newcommand{\C}{\mathbb C} 
\newcommand{\D}{\mathbb D} 
\newcommand{\F}{{\cal F}{}}
\newcommand{\A}{\mathbb A} 
 
\renewcommand{\S}{\mathbb S} 
\newcommand{\M}{\mathbb M}

\newcommand{\R}{\mathbb R}

\newcommand{\N}{\mathbb N}

\newcommand{\sech}{\mbox{sech}} 
 
\newcommand{\arcsinh}{\mbox{\rm arcsinh}} 
\newcommand{\g}{{\mathfrak{g}}{}}

\renewcommand{\k}{{\mathfrak{k}}{}} 
\newcommand{\p}{{\mathfrak{p}}{}} 
 
\newcommand{\n}{{\mathfrak{n}}{}} 
\newcommand{\s}{{\mathfrak{s}}{}}

\renewcommand{\a}{{\mathfrak{a}}{}}

\newcommand{\CO}{{\cal O}{}}
\newcommand{\CM}{{\cal M}{}}

\newcommand{\CS}{\mathcal S}
\newcommand{\CE}{\mathcal E}
\newcommand{\CL}{{\cal L}{}} 
\newcommand{\CD}{\mathcal D}

\newcommand{\CB}{\mathcal B}
\newcommand{\CU}{\mathcal U}
\newcommand{\CF}{\mathcal F}
\newcommand{\CQ}{\mathcal Q}

\def\cref#1{Corollary~\ref{#1}}

\title{Deformation Quantization for actions of 
 the affine group\\
({\sc preliminary version})}

\author{{\bf Pierre Bieliavsky}\\
University of Louvain, Belgium.\\
e-mail: {\tt bieliavsky@math.ucl.ac.be}\\
}
\date{} 

\begin{document}

\maketitle
\begin{abstract}
We define a universal deformation formula (UDF) for the actions of the affine group  on Fr\'echet algebras.
More precisely,
starting with any associative Fr\'echet algebra $\A$ which the affine group $\S\simeq ax+b$ acts on
in a strongly continuous and isometrical manner, the UDF produces a family of topological associative 
algebra structures on the space $\A^\infty$ of smooth vectors of the action deforming the initial
product. The deformation field obtained is based over an infinite dimensional parameter space
naturally associated with the space of pseudo-differential operators on the real line. This note
also presents some geometrical aspects of the UDF and in particular its relation with hyperbolic geometry.
\end{abstract}

\section{Admissible functions on symmetric spaces}
\subsection{Von Neumann's formula in the flat case}
In Weyl's quantization of the flat plane $\R^2$, the formula for the composition
of Weyl's symbols can be expressed in terms of an oscillatory integral three-point kernel
whose phase is proportional to the area of a Euclidean triangle. More precisely,
Weyl's product of two Schwartz symbols $a,b\in \CS(\R^2)$
is given by the following expression, here involving the (real) deformation parameter $\theta$:
\begin{equation}
\left(a\star^{\mbox{\tiny W}}_\theta b\right) (x)=\frac{1}{\theta^2}\int_{\R^2\times\R^2}e^{\frac{i}{\theta}\;S(x,y,z)}\;a(y)\;b(z)\;
{\rm d}y\;{\rm d}z\;.
\end{equation}
In this formula, originally due to Von Neumann \cite{Wi},  $S(x,y,z)$ is a constant multiple of the symplectic area of the Euclidean triangle with vertices $x,y$ and $z$ in $\R^2$. 

\vspace{2mm}

\noindent Only two geometric properties of $S$ are actually sufficient for proving  associtivity of Weyl's product directly at the level of the above formula---i.e. without reference
to any operator representation of symbols \cite{Biep00a}.
The first geometric property is the fact that given any three points $x,y$ and $z$ in the plane, one has:
\begin{equation}\label{ADM1}
S(x,s_x(y),z)=-S(x,y,z)\;,
\end{equation}
where $s_x$ denotes the geodesic symmetry of the plane centred at point $x$.
The second geometric property is the additivity of triangle areas: for all $x,y,z$ and $m$ in the plane, one has:
\begin{equation}\label{COC}
S(x,y,m)+S(y,z,m)+S(z,x,m)=S(x,y,z).
\end{equation}
These properties naturally lead to the following definitions in the more general case of symplectic symmetric spaces.
\subsection{Symplectic  symmetric spaces}
For convenience, we recall in this subsection the notion of symplectic symmetric space ( see \cite{Biep95a,BCG} for details).
\begin{dfn}  A {\bf symplectic symmetric space} is a triple $(M,\omega,s )$ where $(M,\omega )$ is a
connected smooth manifold, $\omega$ is a non-degenerate two-form on $M$ and $s:M\times M\rightarrow M$ is a smooth map such
that \begin{enumerate}	\item[(i)] For any $x$ in $M$, $s_x:M \to M:y \to s(x,y)$ is a $\omega$-preserving
 diffeomorphism of $M$,	which is involutive $(s^2_x=id_{|M})$
and which	admits	$x$ as isolated fixed point. The map $s_x$ is	called	the
{\bf symmetry} at $x$. \item[(ii)] For	any $x,y$ in $M$ one has~:
$$s_xs_ys_x=s_{s_xy}\;.$$ \end{enumerate} 
\end{dfn}
The data of the map $s$ then uniquely determines an affine connection $\nabla$ on $M$ which is invariant under 
the symmetries. The connection turns out to be torsion free and with respect to which the two-form $\omega$ is parallel.
In particular, $\omega$ is symplectic on $M$.

\begin{dfn} 	A	triple	 $t=(\g,\sigma ,{\bf\Omega})$, where $(\g,\sigma)$ is a {\bf involutive Lie algebra}---i.e.
 $\g$ is a finite dimensional real Lie algebra and $\sigma$ is an involutive
automorphism of $\g$ ---and	where	$\bf\Omega$	is	an element of
$\bigwedge^2\g$,	is	called	a {\bf symplectic triple} if the
following properties are satisfied: 
\begin{enumerate}	
\item[(i)]	Let $\g=\k\oplus \p$   where $\k$ (resp. $\p$)
is the $+1$ (resp. $-1$) eigenspace of $\sigma$, then $[\p,\p]=\k$ and the representation of $\k$ on $\p$, given by the adjoint
action, is	faithful. 
\item[(ii)] $\bf\Omega$
is	a Chevalley 2-cocycle for the trivial representation of $\g$ on $\R$
such	that for any	$X$	in	$\k$, $i(X){\bf\Omega}=0$ and such that
the restriction	 of $\bf\Omega$ to $\p \times \p$ defines a symplectic
structure. 
\end{enumerate} 
\end{dfn}

\begin{prop}   There is	a bijection between the set of isomorphism classes of
simply connected symplectic symmetric spaces and the set of isomorphism classes of symplectic triples.
\end{prop}

\noindent Let	us briefly recall how this correspondence is made. If $(M,\omega,s )$ is a
symplectic symmetric space, the group $G$ generated by	products of	an even number of symmetries
is a transitive Lie transformation group of $M$. G is called	the {\bf transvection
group}	of	the symplectic symmetric space.	We associate to	$(M,\omega,s )$ a symplectic triple $t=(\g,\sigma ,{\bf\Omega})$ as follows. The Lie algebra	${\cal	G}$ is	the Lie
algebra	 of	the transvection group	$G$.	Let us choose a	base point $o$ in
$M$ and let $\tilde\sigma$	be	the involutive automorphism of	$G$ obtained
by ``conjugation"	by $s_o$. Then $\sigma$  is the differential of $\tilde\sigma$
at the neutral element	 $e$	of 		$G$. If $\pi :G \rightarrow
M:g\rightarrow g.o$ denotes the projection associated to the choice of base	point
$o$, then $\bf\Omega$	is	the 2-form on $\g$ (identified with $G_e$),
${\bf\Omega}=\pi ^*\omega_o$. The subalgebra $\k$ of $\g$ is known to be isomorphic to
the holonomy algebra for the Loos connection $\nabla$ on the symmetric space
$M=G/_{\textstyle{K}}$---where $K$ is the isotropy at $o$; its Lie algebra is
isomorphic to $\k$. In this context any geodesic through point $o$ is of the form
$t\mapsto\exp(tX)$ where $X\in\p$.

\vspace{3mm}

\noindent In view of UDF's, we will be concerned with symmetric spaces underlying Lie group manifolds \cite{BCSV07}. 
\begin{dfn}
A (symplectic) symmetric space, or more generally a homogeneous space, $M$    
of dimension $m$ is {\bf locally of group type} is there exists a 
$m$-dimensional (symplectic) Lie subgroup $\S$ of its automorphism 
group which acts freely on one of its orbits in $M$.  
One says that it is {\bf globally} of group type if it is locally and    
if $\S$ has only one orbit.
\end{dfn}
In the global case, for every choice 
of a base point $o$ in $M$, the map $\S\to M:g\mapsto g.o$ 
is a $\S$-equivariant diffeomorphism. 
\begin{ex}
{\rm 
\noindent In the present note, we'll be mainly concerned with symmetric {\sl surfaces} (i.e. $\dim M=2$)
 which the affine group $\S=ax+b$ acts on by automorphisms. 
The only globally $\S$-type symmetric surfaces are (see e.g. \cite{Bi1})
the hyperbolic plane $\D:=SL(2,\R)/SO(2)$ and the co-adjoint orbit of the Poincar\'e group $\M:=SO(1,1)\times\R^2/\R$. In particular, one has the following $\S$-equivariant
symplectic identifications: $\M=\S=\D$ and one may therfore think to the space $\M$ as a `curvature contraction'
of the hyperbolic plane $\D$. 
}
\end{ex}
We now detail the above example and denote by $\s$ the (solvable) non-Abelian two-dimensional real Lie algebra which we 
present as generated by the elements $H$ and $E$ with relations $[H,E]=2E$.
Setting $\a:=\R\;H$ and $\n:=\R\;E$ realizes $\s$ as the semi-direct product $\s=\a\times \n$.
The corresponding connected
simply connected Lie group $\S$ is of the solvable exponential type
and the map
\begin{equation}\label{AL}
\s\to\S:(a,\ell):=aH+\ell E\mapsto\exp(aH).\exp(\ell E)
\end{equation}
is a global diffeomorphism.
Within these notations, the group law reads as
\begin{equation*}
(a,\ell).(a',\ell')=(a+a',e^{-2a'}\ell+\ell')\;.
\end{equation*}

\begin{lem}
Consider the solvable symmetric surface $\M=SO(1,1)\times\R^2/\R$. Then,
\begin{enumerate}
\item[(i)] under the identification $\M=\S$, the coordinate system (\ref{AL}) is Darboux ($\omega$ is proportional to ${\rm d}a\wedge{\rm d}\ell$)
and the symmetry map reads:
\begin{equation*}
s_{(a,\ell)}(a',\ell')=\left(2a-a'\;,\;2\cosh(2(a-a'))\ell-\ell'\right).
\end{equation*}
\item[(iii)] The holonomy group $K$ at $o:=(0,0)$ is isomorphic to $\R$ and its
action reads:
\begin{equation*}
\kappa.(a,\ell)=(\;a\;,\;\ell\;-\;\kappa\,\sinh(a)\;)\;.
\end{equation*}
\end{enumerate}
\end{lem}
\subsection{Admissibility}
To a symmetric space $M$, one may attach a natural `group-like' cohomological complex on multiple-point
functions.
\begin{dfn}
Let us define the $k^{\mbox{\rm -th}}$ co-chain space ${\mbox{\rm CP}}^k(M)$ as the space of all complex valued smooth functions on $M^k$ that are invariant under the (diagonal) action of the symmetries on $M^k$. Then,
the formula
\begin{equation*}
\delta F(x_0, ..., x_k)\;:=\;\sum_i(-1)^i F(x_0, ..., \hat{x_i},...,x_k)
\end{equation*}
defines a cohomology operator $\delta: {\mbox{\rm CP}}^k(M)\to {\mbox{\rm CP}}^{k+1}(M)$.
\end{dfn}
Observe that, in this framework, property (\ref{COC}) simply amounts to cocyclicity of $S\in {\mbox{\rm CP}}^3(\R^2)$.
Regarding property (\ref{ADM1}) we observe that, given a geodesic $\gamma$ and given a point $y$ in a symmetric space, 
the curve traced out by $s_xs_z(y)$ where $x$ and $z$ run in $\gamma$ is in general not a geodesic. It is rather the orbit of a one-parameter
transvection subgroup $\{A(t)\}_{t\in\R}$, antifixed under the conjugation by $s_z$, and  realizing the geodesic $\gamma$:
\begin{equation*}
A(t).z=\gamma(t)\;.
\end{equation*}
We then make the following definition, slightly stronger than (\ref{ADM1}).
\begin{dfn} Let $M$ be a geodesically convex symmetric space. A three-cochain $S\in {\mbox{\rm CP}}^3(M)$ is called {\bf admissible} if additionally to $(\ref{ADM1})$ one has,  for all $x,y,z\in M$, that:
\begin{equation}\label{ADM}
S(x, A(t).y,z)=S(x,y,z)\;.
\end{equation}
\end{dfn}
The following fact stresses the relevance of admissible functions on symmetric spaces (for details, see
\cite{Biep00a}). 
\begin{prop} Consider a geodesically convex symplectic symmetric space $M$ and assume
$S\in {\mbox{\rm CP}}^3(M)$ is admissible and cocyclic. Then the product defined by the following formula:
\begin{equation*}
\left(a\star_\theta b\right) (x)=\frac{1}{\theta^2}\int_{M\times M}e^{\frac{i}{\theta}\;S(x,y,z)}\;a(y)\;b(z)\;
{\rm d}y\;{\rm d}z\;\qquad(a,b\in C^\infty_c(M))
\end{equation*}
is formally associative.
\end{prop}
However, 
when curvature is present, admissibility and cocyclicity are, in the non-degenerate case (i.e. 
$S$ of Morse type), incompatible conditions.
Hence, following a standard paradigm in theoretical physics, the idea is to start from an admissible non-cocyclic function and then 
define a deformation framework for it where associativity holds.

\noindent As a starting point in this program, we now give a convenient description of admissible functions on 
strictly geodesically convex symmetric surfaces in terms of admissible functions on the flat plane.

\begin{ex}
{\rm 
Let $(M,\omega,s)$ be a simply connected strictly geodesically convex symplectic symmetric space of dimension two.  
It can then be realized as a coadjoint orbit $\CO$ of a three dimensional Lie group $G$ in the dual 
$\g^\star$ of its Lie algebra $\g$. Denoting by $\g=\k\oplus\p$ 
the decomposition into ($\pm1$)-eigenspaces of $\sigma$, the inclusion $\p\to\g$
induces a canonical projection $\Pi:\g^\star\to\p^\star$ whose restriction to $\CO$ defines, in this case, 
a global
diffeomorphism $$\Pi:\CO\stackrel{\sim}{\to}\p^\star\;.$$
Denote by $o$ the point of $\CO$ corresponding to the origin $0$ of $\p^\star$ under the diffeomorphism 
$\Pi$, and denote by $K$ the stabilizer of $o$ in $G$. For every $x\in\CO$ consider  the associated (globally well-defined) mid-point
map $x\to\frac{x}{2}$ defined by the following property:
$$
s_{\frac{x}{2}}(o):=x\;.
$$
\begin{prop}
View $(\p^\star,\Omega)$ as the flat symplectic plane and consider a $K$-invariant admissible function 
$S^0$ on $\p^\star$ (with respect to the flat structure). Then, the formula:
\begin{equation}\label{ADMFLAT}
S(x,y,z)\;:= \;S^0(0,\Pi(s_{\frac{x}{2}}(y)),\Pi(s_{\frac{x}{2}}(z)))
\end{equation}
defines an admissible function on $M=\CO$. Moreover, every admissible function on $M$
is obtained this way.
\end{prop}
\Pf
We first observe that the 
diffeomorphism $\Pi$ establishes a bijection between the 
$\exp(tX)$-orbits ($X\in\p$) in $\CO$ and the straight lines
in $\p^\star$. Indeed, for $x\in\CO$ and $X\in\p$, one has $<\Ad^\star(\exp(tX))x-x,X>=0$. Which means
that the $x$-translated  $\exp(tX)$-orbit of $x$ lies in the plane in $\g^\star$ orthodual to $X\in\p$.
This plane is generated by the kernel $\k^\star$ of the projection $\Pi:\g^\star\to\p^\star$ and 
an element $X^\perp$ of $\p^\star$ orthodual to $X$. In particular, it projects onto the line 
directed by $X^\perp$.

\noindent Now consider the two-point function $u(x,y):=S(o,x,y)$ on $M$ induced by the data of an 
admissible function $S$ on $M$. This function corresponds to a two-point function $u^0$ on $\p^\star$
via the diffeomorphism $\Pi$. By admissibility (cf. (\ref{ADM})) and the above observation, one has
$u^0(\xi,\eta)=u^0(\xi,\eta+t\xi)$ for all $t\in\R$. Which is precisely the property of admissibility
for a two-point function with respect to the flat structure on $\p^\star$. The rest then follows from Proposition 3.3 in \cite{Biep00a}.
\EPf
\begin{dfn}
By virtue of the above proposition, the  Euclidean triangle symplectic area on $\p^\star$ defines an admissible
three-point function on $M$. The latter will be denoted $S_{{\rm can}}$ and called the {\sl canonical} admissible function.
\end{dfn}
Observe furthermore that any odd function of a multiple of $S_{\rm can}$ defines an admissible
three-point function on $M$. These particular admissible functions correspond to the ones which in addition
to their $\k$-symmetry
enjoy a full ${\mathfrak {sp}}(1,\R)$-symmetry.
}
\end{ex}
\begin{ex}
Within the coordinate system (\ref{AL}), the canonical admissible function $S_{{\rm can}}$ on $\M=SO(1,1)\times\R^2/\R$
has the following expression:
\begin{equation}
S_{{\rm can}}(x_0,x_1,x_2)=\frac{1}{2}\cyclic_{0,1,2}\sinh(2(a_0-a_1))\ell_2\;,
\end{equation}
where $x_i=(a_i,\ell_i)\quad (i=0,1,2)$.
\end{ex}
\section{Strict quantizations of $SO(1,1)\times\R^2/\R$}\label{SOSS}
In this section, we describe all the invariant deformation quantizations on the symplectic symmetric surface $\M:=SO(1,1)\times\R^2/\R$. We begin by slightly generalizing the WKB-quantization constructed in \cite{Biep00a}.
What follows also provides a rigourous framework for statements made in \cite{BDRS}.

\noindent Endowing 
$\S$ with any left-invariant Haar measure, the Darboux map (\ref{AL}) yields the identifications
\begin{equation*}
L^2(\S)=L^2(\s)=L^2(\M)\;.
\end{equation*}
The partial Fourier transform
with respect to the  variable $\ell$ will be denoted by
\begin{eqnarray*}
&\F_N:L^2(\S)\to L^2(\tilde{\S})& \\
&\F_N(u)(a,\alpha):=\int_\n e^{-i\alpha \ell}u(a,\ell)\;{\rm d}\ell\;,&
\end{eqnarray*}
where $\tilde{\S}:=\{(a,\alpha)\}$ denotes the space where Fourier transformed functions are defined on.
\begin{dfn}
The {\bf twisting map} is the one-parameter family of diffeomorphisms of $\tilde{\S}$ defined by
\begin{equation*}
\varphi_\theta(a,\alpha):=(\;a\;,\;\frac{1}{2\theta}\sinh\left(2\theta\alpha\right)\;)\qquad \theta\in\R.
\end{equation*}
\end{dfn}
Observe that, denoting by $\tilde{\CS}$ the Schwartz function space on $\tilde{\S}$, 
on has 
\begin{equation*}
\varphi_\theta^\star(\tilde{\CS})\subset\tilde{\CS}\;.
\end{equation*}
Therefore, for every invertible operator multiplier $\Theta\in\CO_M(\R)$, the following linear map defined on the Schwartz space $\CS$ on $\S$ (with respect to coordinates $(a,\ell)$) takes its values 
in the tempered  smooth functions on $\S$:
\begin{equation}\label{EQ0}
U^\Theta_\theta\;:=\;\CF_N^{-1}\circ\CM_\Theta\circ(\varphi_\theta^{-1})^\star\circ\CF_N\;:\;\CS\;\to\;\CS'\;.
\end{equation}
Its range will be denoted
\begin{equation*}
\CE^\Theta_\theta:=U^\Theta_\theta(\CS)\subset\CS'\;;
\end{equation*}
and we set  
\begin{equation*}
\left(U^\Theta_\theta\right)^{-1}:\CE^\Theta_\theta\stackrel{\sim}{\to}\CS
\end{equation*}
for the associated (inverse) linear isomorphism. The space $\CE^\Theta_\theta$ carries
the transported Schwartz topology. Note that following the same argument as in \cite{Biep00a}, one has the inclusion:
\begin{equation*}
\CE^\Theta_\theta\supset\CS\;.
\end{equation*}
\begin{thm}\label{CONTRACTED}
Let $\Theta$ be any invertible element of $\CO_M(\R)$ and consider $\theta>0$.
\begin{enumerate}
\item[(i)] Let $u$ and $v$ in $C^\infty_c(\M)$ and $x_{0}\in\M$. Then, the formula\footnote{Formula (i) below has been announced in \cite{BDRS}.}:
\begin{equation*}
(u{\star}^\Theta_\theta \,v)(x_0)\;:=\;\frac{1}{2\pi\theta^2}\int_{\M\times \M}{{ \mbox{$\cosh(2(a_1-a_2))
\frac{{\bf \Xi}(a_2-a_0){\bf \Xi}(a_0-a_1)}{{\bf \Xi}(a_2-a_1)}$}}}\;e^{\frac{i}{\theta}S_{{\rm can}}(x_0,x_1,x_2)}\;u(x_1)\;v(x_2)\;{\rm d}x_1\;{\rm d}x_2\;,
\end{equation*}
where 
\begin{equation*}
{\bf \Xi}(\theta t)\;:=\;\frac{1}{(\varphi_\theta^\star\Theta)(t)}\;,
\end{equation*}
extends to $\CE^\Theta_\theta$ as an associative product. The pair $(\CE^\Theta_\theta,{\star}^\Theta_\theta)$
is then a Fr\'echet algebra.
\item[(ii)]  The formula
\begin{equation*}
\mbox{Tr}^\Theta_{\theta}¥(u)\;=\;\frac{1}{\Theta(0)}\int_\M u\;\qquad (u\in\CS)
\end{equation*}
extends to $\CE^\Theta_\theta$ as a trace
$\mbox{Tr}^\Theta_{\theta}¥:\CE^\Theta_{\theta}¥\to\C$ for ${\star}^\Theta_\theta$.
\item[(iii)] The algebra $(\CE_{\theta}^{\Theta},¥¥{\star}^\Theta_\theta)$ is strongly closed if and only if
\begin{equation}\label{UNIT}
{\bf \Xi}(t){\bf \Xi}^{\vee}(-t)\;=\;\cosh(2t)\;.
\end{equation}
\end{enumerate}
\end{thm}
\Pf
Item (i) is obtained by a long but straightforward computation, entirely similar to the one in
\cite{Biep00a} (see the proof of Theorem 6.13 page 311), where one transports Weyl's 
product on $\CS=\CS_{(a,\ell)}$ to $\CE^\Theta_{\theta}$ via $U^\Theta_\theta$.

\noindent The trace formula is obtained by combining 
$\mbox{Tr}^\Theta=\mbox{Tr}^{\mbox{{\tiny{W}}}}\circ\left(U^\Theta_\theta\right)^{-1}$ where $\mbox{Tr}^{\mbox{{\tiny{W}}}}(u):=\int_\M u$ with the fact that $\int\circ\,\CF^{-1}_N=\delta_0$.
At last, using strong closedness of Weyl's product, one gets
\begin{eqnarray*}
\mbox{Tr}^\Theta(u{\star}^\Theta_\theta \,v)=\int\left(U^\Theta_\theta\right)^{-1}(u)\;\left(U^\Theta_\theta\right)^{-1}(v).
\end{eqnarray*}
Which equals $<\left(U^\Theta_\theta\right)^{-1}(u)\,,\,\overline{\left(U^\Theta_\theta\right)^{-1}(v)}>_{L^2(\M)}$.
Setting $f_{\theta}(t):=f(\theta t)$ and $f^\vee(t):=f(-t)$, the latter becomes
\begin{eqnarray*}
<\CF_N^{-1}\left(\Theta_\theta\,\varphi_\theta^\star\CF_Nu\right)\,,
\,\CF_N^{-1}\left[\overline{\left(\Theta_\theta\,\varphi_\theta^\star\CF_Nu\right)}\right]^{\vee}>_{L^2(\M)}=
\int\Theta_\theta\Theta_\theta^\vee\,\varphi_\theta^\star\left(\CF_N(u)\,\left[{\CF_N(v)}\right]^\vee\right)=\\
=\int(\varphi^{-1}_\theta)^\star(\Theta_\theta\Theta_\theta^\vee)\,.\,|\mbox{\rm Jac}_{\varphi^{-1}_\theta}|\,.\,\CF_N(u)\,\left[{\CF_N(v)}\right]^\vee.
\end{eqnarray*}
The last member equals $\int uv$ if and only if 
$(\varphi^{-1}_\theta)^\star(\Theta_\theta\Theta_\theta^\vee)\,.\,|\mbox{\rm Jac}_{\varphi^{-1}_\theta}|=1$, which is the announced condition.\EPf
\begin{rmk}
{\rm 
\noindent (i)  The invertible element $\Theta$ can be considered as a parameter in the construction. Moreover,
it can itself depend on the real parameter $\theta$ as well.

\noindent (ii) By defining $<a,b>_\Theta:=\mbox{Tr}^\Theta(a\star^\Theta\overline{b})$, one could study
the associated field of Hilbert $G(\M)$-algebras in the same line as in \cite{Biep00a}. We will not follow this route
here, but rather focus on the square integrable case associated with the unitary condition (\ref{UNIT}).
}
\end{rmk}
\noindent The manifold $\M\times \M\times \M$ admits a distinguished transformation, namely one has
\begin{lem}\cite{Qi97,Biep00a}
\begin{enumerate}
\item[(i)] Given any triple of points $x,y$ and $z$ in $\M$, the equation
\begin{equation*}
s_xs_ys_z(t)=t
\end{equation*}
admits a unique solution $t\in \M$.
\item[(ii)] The associated map
\begin{equation*}
\Phi:\M\times \M\times \M\to \M\times \M\times \M:(x,y,z)\mapsto(t,s_z(t),s_ys_z(t))
\end{equation*}
is a global diffeomorphism.
\end{enumerate}
\end{lem}
\begin{lem}
\begin{equation*}
\mbox{\rm Jac}_{\Phi}(x_{0},x_1,x_2)\;=\;16\,\cosh(2(a_0-a_1))\,\cosh(2(a_1-a_2))\,\cosh(2(a_2-a_0))\;.
\end{equation*}
\end{lem}
\Pf
This is a straightforward computation based on the following formula for the mid-point map:
\begin{eqnarray*}
m:M\times M\to M:(x,y)\mapsto m(x,y)=\left(\;\frac{1}{2}(a_x+a_y)\;,\;\frac{1}{2}(\ell_x+\ell_y)\sech(2(a_x-a_y))\;\right)\;,
\end{eqnarray*}
defined by the relation
\begin{eqnarray*}
s_{m(x,y)}x=y\;.
\end{eqnarray*}
One has 
\begin{equation*}
\Phi^{-1}(x,y,z)=(m(x,y),m(y,z),m(z,x))\;;
\end{equation*}
and a computation yields
\begin{equation*}
\mbox{\rm Jac}_{\Phi^{-1}}(x_0,x_1,x_2)=\frac{1}{16}\sech(2(a_0-a_1))\sech(2(a_1-a_2))\sech(2(a_2-a_0)).
\end{equation*}
One then obtains the announced formula by using the relation: 
$\mbox{\rm Jac}_{\Phi}=(\Phi^\star\mbox{\rm Jac}_{\Phi^{-1}})^{-1}$.
\EPf
The latter together with the unitary condition (\ref{UNIT}) yield
\begin{cor}\label{THM1}
Let $\theta>0$.  Then,
\begin{enumerate}
\item[(i)]  the following formula
\begin{equation}\label{star1} 
u\star_\theta v(x):=\frac{1}{\theta^2}\int_{\M\times \M}\left[\mbox{Jac}_\Phi(x,y,z)\right]^\frac{1}{2}
\;e^{\frac{i}{\theta}S_{\rm can}(x,y,z)}\;u(y)\;v(z)\;{\rm d}y\;{\rm d}z\qquad(u,v\in C^\infty_c(\M))
\end{equation}
extends to the space\footnote{$L^2(\M)$ denotes the space of square-integrable functions with respect to the Liouville measure.} 
$L^2(\M)$ as an associative product.
\item[(ii)] The algebra $(L^2(\M),\star_\theta)$ becomes a Hilbert algebra when one endows
$L^2(\M)$ with its natural Hilbert space structure.
\item[(iii)] The transvection group $G=G(\M)$ acts on the above algebra by unitary automorphisms.
\end{enumerate}
\end{cor}

\noindent remains to give a geometrical meaning to the co-boundary factor appearing in the prpoduct formula. 
For this, we observe that the affine manifold $\M$ admits a canonical invariant foliation with respect to which $\ell$-independent 
functions on $\M$  correspond to leafwise constant functions.
\begin{lem}\label{FOL}\cite{Biep00a}
The symmetric surface $\M$ admits a unique one-dimensional distribution $\CL\subset T(\M)$ which is invariant
under the transvection group. The corresponding (Lagrangian) foliation of $\M$ is a fibration by geodesics.
\end{lem}
The twisting map now appears as a one parameter family of transformations of $\M/\CL$. Moreover, one notes
\begin{cor}
The  factor $\frac{{\bf \Xi}(a_2-a_0){\bf \Xi}(a_0-a_1)}{{\bf \Xi}(a_2-a_1)}$ corresponds
to  $ \exp(\delta\Xi)$ where $\Xi$ is a symmetry invariant $\CL\times\CL$-leafwise constant 2-cochain in ${\mbox{\rm CP}^2(\M)}$.
\end{cor}
We end this section with a remark concerning invariant  formal star products on $\M$. First we observe that following  the same lines as 
in the proof of Proposition 4.2. in \cite{BBM}, one obtains
\begin{prop}
Consider a $\CO_M(\R)$-valued smooth function $\Theta:]-\epsilon,\epsilon[\to\CO_M(\R):\theta\mapsto\Theta_\theta$ such that
$\Theta_\theta$ is invertible for all $\theta$ and such that $\Theta_0\equiv1$. Let $u$ and $v$ be smooth and compactly supported and consider the associated 
product 
\begin{equation*}
u\star^\Theta_\theta v=\frac{1}{\theta^2}\int_{\M\times\M}\,|{\mbox{\rm Jac}}_\Phi|^{\frac{1}{2}}\,\exp(\delta\Xi_\theta)\,e^{\frac{i}{\theta}S_{\mbox{\rm can}}}\,u\otimes v\;,
\end{equation*}
 as in Theorem \ref{CONTRACTED} (i). Then the function $\theta\mapsto u\star^\Theta_\theta v$ is a smooth $C^\infty(\M)$-valued function whose Taylor series at $\theta=0$
defines a symmetry invariant  formal star product $\tilde{\star}^\Theta_\theta$ on $C^\infty(\M)[[\theta]]$.
\end{prop}
The remark is then  that every invariant star product on $\M$ can be seen as an asymptotic expansion
of an oscillatory integral of the above type. Indeed,
through the identification
$\S=M$, the formal star product $\tilde{\star}^1_\theta$  can be viewed as
a left-invariant formal star product on the symplectic Lie group $\S$ with an additional $K$-invariance.
Every other $G$-invariant star product on $(\M,\omega)$ in the same $G$-characteristic class can therefore be obtained by intertwining
$\tilde{\star}^1_\theta$ by a left-invariant formal equivalence (i.e. an element of $\CU(\s)[[\theta]]$) which commutes with the $K$-action.
Within the coordinate system (\ref{AL}), the left-invariant vector fields on $\S$ corresponding to elements
$H$ and $E$ of $\s$ have the following expressions:
\begin{equation*}
\widetilde{H}=\partial_a-2\ell\partial_\ell\;,\mbox{ and }\;\widetilde{E}=\partial_\ell.
\end{equation*}
While the fundamental vector field
\begin{equation*}
Z_c^\star:=\sinh(a)\partial_\ell
\end{equation*}
is a generator of the $K$-action.
From this, one sees that each term of the above mentioned equivalence must be  polynomial in $\widetilde{E}$ only. In other words, 
the conjugation by $\CF_N$ of the equivalence is a multiplication by an (inversible)
formal function. Combining the above observation with Formula (\ref{EQ0}) basically yields the following proposition \begin{prop}
Let $\Theta\in \C[t][[\theta]]$ with non-zero constant leading term. Then, 
the formal star product defined as 
\begin{equation*}
\tilde{\star}^\Theta_\theta := \CF_N^{-1}\circ\CM_{\Theta}\circ\CF_N(\tilde{\star}^1_\theta)
\end{equation*}
is  $G$-invariant (in the above formula $\CM_\Theta$ denotes 
the multiplication operator: $\CM_\Theta(f):=\Theta f$). Moreover, every $G$-invariant
formal star product on $(\M,\omega)$ is of this form.
\end{prop}
\Pf
It remains to prove the second assertion, which follows from the fact that changing $G$-characteristic class
amounts to change the parameter.  Indeed, one knows that the $G$-equivalence classes of 
invariant star products on $(\M,\omega)$ are canonically parametrized by the formal series
with coefficients in the second $G$-equivariant de Rham cohomology space $H^2_{G}(\M)$ \cite{BBG}. The latter
is isomorphic to $\C$ as generated by $[\omega]$. Hence every $G$-invariant class is of the 
form $\Theta([\omega])$ where $\Theta=\Theta(\theta)$ is a formal function. One then concludes by
Proposition 4.3 in \cite{BiBo}.
\EPf
\section{Universal deformation Formulae}
\subsection{Oscillatory integrals} We let $E$ be a complex Fr\'echet space with topology defining family of seminorms $\{|\;\;|_j\}_{j\in\N}$.
Let $G$ be a solvable exponential Lie group with Lie algebra $\g$
and consider a function ${\bf m} \in C^\infty(G,\C)$ which is nowhere vanishing.
We then define the following function space:
\begin{equation*}
\CB^{\bf m}_E(G) := \{F\in C^\infty(G, E)\mbox{ such that }\forall P\in\CU(\g)\;;\;j\in\N\;\mbox{ there exists }C>0\mbox{ such that  }\;|\tilde{P}.F|_j<C\,|{\bf m}|\,\}\;;
\end{equation*} 
where $\tilde{P}$ denotes the left-invariant differential operator on $G$ associated with the element $P$
of the universal enveloping algebra $\CU(\g)$.
\begin{dfn}
An everywhere non-zero  function ${\bf m}\in C^\infty(G,\C)$ is called a {\bf weight} if ${\bf m}\in\CB^{\bf m}_\C$.
\end{dfn}
Let $C_b(G,E)$ be the space of $E$-valued bounded continuous functions on $G$. The group
$G$ then acts on the latter space  via the right regular representation. Consider the subspace
$C_u(G, E)$ of $C_b(G,E)$  constituted by the uniformly continuous functions. The following lemma is essentially
standard.
\begin{lem}\label{symbols}
Let ${\bf m} $ and ${\bf m}'$ be weight functions. Then,
\begin{enumerate}
\item[(i)] the group $G$ acts on $C_u(G, E)$ isometrically and strongly continuously. The space $[C_u(G, E)]_\infty$
of smooth vectors for this action coincides with  $\CB^1_E(G)$.
\item[(ii)] On $\CB^1_E(G)$,  the following seminorms:
\begin{equation*}
|a|_{P,j}\;:=\;\sup_G\{|\tilde{P}.a|_j\}\qquad(P\in\CU(\g)),
\end{equation*}
induce the natural Fr\'echet topology on $[C_u(G, E)]_\infty$ (cf. \cite{War}).

\noindent Analogously, on $\CB^{\bf m}_E(G)$ the seminorms
\begin{equation*}
|a|_{P,j}\;:=\;\sup_G\{|\frac{1}{\bf m}\,\tilde{P}.a|_j\}
\end{equation*}
define a Fr\'echet topology.

\item[(iii)] The group $G$ acts isometrically on $\CB^1_E$ via the left regular representation. In particular, the space $[C_u(G, E)]_\infty$
is a $G$-bimodule.
\item[(iv)] For every $u\in\CB^{{\bf m}}_\C$ and $a\in\CB^{{\bf m}'}_E$, their product, $ua$, belongs to $\CB^{{\bf m}{\bf m}'}_E$. Moreover,
the associated bilinear map:
\begin{equation*}
\CB^{{\bf m}}_\C\times\CB^{{\bf m}'}_E\to\CB^{{\bf m}{\bf m}'}_E
\end{equation*}
is continuous.
\item[(v)] For every $P\in\CU(\g)$ and $a\in\CB^{\bf m}_E$, the element $\tilde{P}.a$ belongs to $\CB^{\bf m}_E$ and the 
map
\begin{equation*}
\CB^{\bf m}_E\to\CB^{\bf m}_E: a\mapsto\tilde{P}.a
\end{equation*}
is continuous.
\item[(vi)] Assume the inverse, $\frac{1}{\bf m}$,  of the weight  ${\bf m}$ vanishes at infinity. Then, the closure of $\CD_E$
in $\CB^{\bf m}_E$ contains $\CB^{{\bf m}'}_E$ for all ${\bf m}'$ such that $|{\bf m}'|<|{\bf m}|$.
\item[(vii)] 
Let ${\bf m}\in\CB^{\bf m}_\C$ be a weight and consider any nowhere vanishing function ${\bf m}_0$. Then for all $A\in\CB^{{\bf m}_0}_E$ and all $P\in\CU(\g)$,
one has $\tilde{P}.(\frac{A}{\bf m})=\frac{1}{\bf m}\,A'$ where $A'$ belongs to $\CB^{{\bf m}_0}_E$. 
\end{enumerate}
\end{lem}
\Pf
 The space
$C_u(G, E)$ is a Fr\'echet space for the seminorms $\{|\;|_j^\infty\}$ defined as 
\begin{equation*}
|a|^\infty_j\;:=\;\sup_G|a|_j\;.
\end{equation*}
Indeed, $G$ being locally compact and countable at infinity the space $C_b(G,E)$ is Fr\'echet (by the same argument
as in the proof of Prop. 44.1 and Cor. 1. of \cite{Tr}). The subspace $C_u(G, E)$ is then closed as 
a uniform limit of uniformly continuous functions is uniformly continuous (as it is seen by a $3$-epsilon 
argument). 

\noindent  Moreover, the natural Fr\'echet topology 
on $[C_u(G, E)]_\infty$ is induced by the set of seminorms $\{|\;|_{P,j}\}$ on $\CB^1_E(G)$ defined as 
\begin{equation*}
|a|_{P,j}\;:=\;\sup_{G}|\tilde{P}.F|_j\;.
\end{equation*}
Indeed, an element $a\in [C_u(G, E)]_\infty$ is such that the function $g\mapsto R^\star_ga$ is smooth as a $C_u(G, E)$-valued function
on $G$. In particular, for every $P\in\CU(\g)$, $\tilde{P}.a$ is bounded and smooth. Reciprocally,
$G$ acts on $\CB^1_E$ via the right regular representation. Indeed,  for all $g\in G, \tilde{P}_x(R^\star_ga)=$ $(\Ad(g^{-1})P)^\sim|_{xg}.a$.
Hence $\sup_x|\tilde{P}_x(R^\star_ga)|_j=\sup_{x}|(\Ad(g^{-1})P)^\sim|_{xg}.a|_j=\sup_{x}|(\Ad(g^{-1})P)^\sim|_{x}.a|_j$ which is bounded
for $a\in\CB^1_E$. Note that the group $G$ acts on $\CB^1_E$ via the left regular representation as well. Indeed, $\sup|\tilde{P}.(L_g^\star a)|_j=
\sup|L_g^\star(\tilde{P}.a)|_j=\sup|\tilde{P}.a|_j$. The left regular representation action is in particular isometric.
Note also that one has the inclusion:
$\CB^1_E\subset C_u(G,E)$. Indeed,  for $a\in\CB^1_\C$ the function $ \tilde{{\rm d} a}:G\to\g^\star:x\mapsto{\rm d}a_x\circ L_{x\star_e}$ is such that
$<\tilde{{\rm d} a}(x),H>\leq c(H)$ where $c:\g\to \C$ is independent of $x$. One may moreover assume that the function $c$
is continuous. Indeed, setting $c_s(H):=\sup_x\{|<\tilde{{\rm d} a}(x),H>|\}$, one observes that $c_s(\lambda H)=|\lambda| c_s(H)$ and 
$c_s(H+H')\leq c_s(H)+c_s(H')$. Choosing a basis $\{X_j\}$ of $\g$, one then gets positive numbers $\{m_j\}$ such that $c_s(x^jX_j)\leq\sum_jm_j\,|x^j|\,=:c(x^jX_j)$.
Now, for fixed $H$, one observes that $|a(x\exp(tH))-a(x)|=|{\rm d}a_{x\exp(\tau H)}(\tilde{H})\,t|=|<\tilde{{\rm d}a}(x\exp(\tau H)),H>t|$ where $\tau\in[0,t]$. Therefore
$|a(x\exp(tH))-a(x)|\leq c(H)|t|$. Choosing a Euclidean scalar product on $\g$, and denoting by $B_r$ the open ball of radius $r$ in $\g$, one 
observes  that for all $x$ in $G$, one has
$|a(x\exp(B_r))-a(x)|\leq\max_{|H|=1}(c)\,r$; hence the uniform continuity of $a$. To show that $a\in\CB^1_\C$ is a differentiable vector, we observe that
\begin{eqnarray*}
\sup_x\{|\frac{1}{t}(a(x\exp(t X))-a(x))-\tilde{X}_x.a|\}&=&\sup_x\{|\tilde{X}_{x\exp(\tau X)}.a-\tilde{X}_x.a|\} \qquad(\tau\in[0,t])\\
 &=&|\tau|\,\sup_x\{|\tilde{X}^2_{x\exp(\sigma X)}.a|\}\qquad(\sigma\in[0,\tau])\\  &\leq&|t|\sup\{|\tilde{X}^2.a|\}\,,
\end{eqnarray*}
which tends to zero together with $t$. This yields differentiability at the unit element. One gets it everywhere else by observing that
\begin{equation}\label{tech1}
\tilde{X}. R^\star(a)=R^\star(\tilde{X}.a)\;.
\end{equation}
An induction on the order of derivation implies $\CB^1_\C\subset [C_u(G)]_\infty$,
and the $E$-valued case is entirely similar. 

\noindent The assertion concerning the topology follows from the definition of the topology on smooth vectors \cite{War} and from
(\ref{tech1}) again. 

\noindent Now, the notion of weight  implies that $a$ belongs to $\CB^{\bf m}_E$ iff $\frac{1}{\bf m}a$ belongs to $\CB^1_E$; the
non-constant weight case then follows.

\noindent Items (iv) and (vii) follows from Leibniz' rule while item (v) is obvious.

\noindent For the last assertion, we consider, similarly as in \cite{Di}, a cut-off $\varphi\in\CD$
such that $\varphi|_{B_1}=1$ and $\varphi|_{G\backslash B_2}=0$. Setting $\varphi_n(x):=\varphi(\frac{1}{n}x)$,
we observe that $\sup\{|\frac{1}{\bf m}\,\tilde{P}.(1-\varphi_n)|_j\}$ tends to zero when $n$ tends to infinity for every
$P\in\CU(\g)$. Which amounts to say that $\{\varphi_n\}$ converges to $1$ in $\CB^{\bf m}_\C$.
The latter combined with item (iv) yield item (vi). \EPf
\begin{dfn} Consider  a function $S\in C^\infty(G,\R)$.
An element ${\bf P}$ of $\CU(\g)$ is called {\bf $S$-adapted} if the following conditions hold
\begin{enumerate}
\item[(i)] the function ${\bf m}_{\bf P}:= e^{-iS}\tilde{\bf P}.e^{iS}$ is a  weight;
\item[(ii)] its inverse $\frac{1}{{\bf m}_{\bf P}}$ is integrable with respect to a left-invariant Haar measure on $G$.
\end{enumerate}
If such an element ${\bf P}$ exists then one calls $S$ a {\bf phase} on $G$. Moreover, for every weight ${\bf m}_0$ such that $\frac{{\bf m}_0}{{\bf m}_{\bf P}}$ is integrable,
we call an {\bf amplitude} (adapted to $S$) any element of $\CB^{{\bf m}_0}_E$.
\end{dfn}
\begin{rmk}
{\rm 
Note that the product of an amplitude by an element of $\CB^1_\C$ is again an amplitude.
}
\end{rmk}
One then has
\begin{dfn} Assume ${\bf P}\in\CU(\g)$ is $S$-adapted and self-adjoint\footnote{The latter condition of self-adjointness stands there
for simplicity, but it is not essential.}.
Consider any  weight ${\bf m}_0$ such that $\frac{{\bf m}_0}{{\bf m}_{\bf P}}$ is integrable. 
For $A\in\CD_E$, an integration\footnote{For standard definitions and properties about $L^1(G,E)$, we refer
to \cite{Gr}.} by parts yields $\int_G e^{iS}A=\int_Ge^{iS}\tilde{\bf P}.(\frac{A}{{\bf m}_{\bf P}})$. Moreover,
Lemma \ref{symbols}  implies that
the linear map $\CD_E\to E:A\mapsto\int e^{iS}A$ extends by continuity to a continuous linear map:
\begin{equation}\label{osc}
\tilde{\int}\,e^{iS}\,:\CB^{{\bf m}_0}_E\to E\;.
\end{equation}
The latter is called the {\bf oscillatory integral} on $\CB^{{\bf m}_0}_E$.
\end{dfn}
\subsection{Hyperbolic Laplacian and the canonical phase}
\begin{lem}
On the affine group $\S$, one has 
\begin{equation*}
{}^\tau\tilde{E}=-\tilde{E}\;;{}^\tau\tilde{H}=-\tilde{H}+2\;;
\end{equation*}
and the operator
\begin{equation*}
\tilde{B}:=\alpha\tilde{H}^2+\beta\tilde{E}^2+\gamma\frac{1}{2}(\tilde{E}\tilde{H}+\tilde{H}\tilde{E})-2\alpha\tilde{H}
-\gamma\tilde{E}
\end{equation*}
is self-adjoint for every data of  $\alpha,\beta,\gamma$ in $\R$.
\end{lem}
In particular, the Laplace operator associated with any  Damek-Ricci Riemannian structure on $\S$ (\cite{DaRi}) is
of the above form. The analysis which follows may be performed with any of the latter, however, for simplicity we shall only consider the hyperbolic Laplacian:
\begin{equation*}
\Delta\;:=\;2(\,\tilde{H}^2+\tilde{E}^2-2\tilde{H}\,)\;.
\end{equation*}
\begin{prop}\label{laplacian} The canonical  two-point function on $\S=M$:
\begin{equation*}
S_{\mbox{\rm can}}\;:=\;\sinh(2a_1)\ell_2-\sinh(2a_2)\ell_1\;
\end{equation*}
is a phase on $\S\times\S$. 

\noindent Moreover, every weight ${\bf m}_0={\bf m}_0(a_1,a_2)$ such that 
\begin{equation*}
\frac{{\bf m}_0}{\sqrt{\cosh(2(a_1-a_2))\cosh(4(a_1+a_2))\cosh(4(a_1-a_2))}}\in L^1_{(a_1,a_2)}
\end{equation*}
is an adapted amplitude.
\end{prop} 
\Pf  Define, for $\phi\in C^\infty(\S\times\S)$, 
\begin{equation*}
{\bf \Delta}\,\phi\;:=\;(\,\Delta\otimes1+1\otimes\Delta\,)\,\phi\;.
\end{equation*}
A computation shows that
\begin{equation*}
{\bf \Delta}^2\,e^{i S_{\mbox{\rm can}}}\;=\left(\CQ_4+\CQ_2 +{\bf c}+i\CQ_3\right)\;e^{iS_{\mbox{\rm can}}}\;.
\end{equation*}
where 
\begin{eqnarray*}
\CQ_4\;:=\;128 \left((\begin{array}{cc}
\ell_1&\ell_2
\end{array})
\left(
\begin{array}{cc}
 \cosh[4 a_2] & \sinh[2 (a_1 + a_2)]\\
 \sinh[2 (a_1 + a_2)]& \cosh[4 a_1] 
\end{array}
\right)
\left(
\begin{array}{c}
\ell_1\\\ell_2
\end{array}
\right)
\right)^2\;;
\end{eqnarray*}

\begin{eqnarray*}
\CQ_2\;:=\; 
{\ell_2}^2\,\big(528 - 2752\,\cosh (4\,a_1) + 16\,\cosh (8\,a_1) +  
     16\,\cosh (4\,\left( a_1 - a_2 \right) ) + \\
     16\,\cosh (4\,\left( a_1 + a_2 \right) ) + 896\,\sinh (4\,a_1) \big)  + 
  {\ell_1}^2\,\big( 528 + 16\,\cosh (4\,\left( a_1 - a_2 \right) ) -\\
     2752\,\cosh (4\,a_2) + 16\,\cosh (8\,a_2) + 
     16\,\cosh (4\,\left( a_1 + a_2 \right) ) + 896\,\sinh (4\,a_2) \big)  + \\
  \ell_1\,\ell_2\,\big( -1024\,\cosh (2\,\left( a_1 - a_2 \right) ) + 
     3328\,\cosh (2\,\left( a_1 + a_2 \right) ) - \\
     3968\,\sinh (2\,\left( a_1 + a_2 \right) ) + 
     32\,\sinh (2\,\left( 3\,a_1 + a_2 \right) ) + \\
     32\,\sinh (2\,\left( a_1 + 3\,a_2 \right) ) \big)\;;
\end{eqnarray*}

\begin{eqnarray*}
\CQ_3\;:=\; {\ell_2}^3\,\big( \left( 128\,  \right) \,\cosh (2\,a_1) + \left( 128\,  \right) \,\cosh (6\,a_1) +\\
 \left( 384\,  \right) \,\sinh (2\,a_1) - \left( 896\,  \right) \,\sinh (6\,a_1) \big)  
\\ + {\ell_1}^3\,\big( \left( -128\,  \right) \,\cosh (2\,a_2) - \\
     \left( 128\,  \right) \,\cosh (6\,a_2) - 
     \left( 384\,  \right) \,\sinh (2\,a_2) + \left( 896\,  \right) \,\sinh (6\,a_2)
     \big)  +\\
 \ell_1\,\big( \left( 16\,  \right) \,\cosh (4\,a_1 - 2\,a_2) + 
     \left( 880\,  \right) \,\cosh (2\,a_2) - \\
     \left( 16\,  \right) \,\cosh (6\,a_2) - 
     \left( 112\,  \right) \,\cosh (2\,\left( 2\,a_1 + a_2 \right) ) - \\
     \left( 48\,  \right) \,\sinh (4\,a_1 - 2\,a_2) - 
     \left( 1552\,  \right) \,\sinh (2\,a_2) + \\
     \left( 144\,  \right) \,\sinh (6\,a_2) + 
     \left( 48\,  \right) \,\sinh (2\,\left( 2\,a_1 + a_2 \right) ) \big)  + \\
  \ell_1\,{\ell_2}^2\,\big( \left( -128\,  \right) \,
      \cosh (4\,a_1 - 2\,a_2) + \left( 768\,  \right) \,\cosh (2\,a_2) - \\
     \left( 1920\,  \right) \,\cosh (2\,\left( 2\,a_1 + a_2 \right) ) - 
     \left( 384\,  \right) \,\sinh (4\,a_1 - 2\,a_2) + \\
     \left( 256\,  \right) \,\sinh (2\,a_2) + 
     \left( 1152\,  \right) \,\sinh (2\,\left( 2\,a_1 + a_2 \right) ) \big)  + \\
  \ell_2\,\big( \left( -880\,  \right) \,\cosh (2\,a_1) + 
     \left( 16\,  \right) \,\cosh (6\,a_1) - \\
     \left( 16\,  \right) \,\cosh (2\,\left( a_1 - 2\,a_2 \right) ) + 
     \left( 112\,  \right) \,\cosh (2\,\left( a_1 + 2\,a_2 \right) ) + \\
     \left( 1552\,  \right) \,\sinh (2\,a_1) - 
     \left( 144\,  \right) \,\sinh (6\,a_1) - 
     \left( 48\,  \right) \,\sinh (2\,\left( a_1 - 2\,a_2 \right) ) - \\
     \left( 48\,  \right) \,\sinh (2\,\left( a_1 + 2\,a_2 \right) ) \big)  + 
  {\ell_1}^2\,\ell_2\,\big( \left( -768\,  \right) \,\cosh (2\,a_1) + \\
     \left( 128\,  \right) \,\cosh (2\,\left( a_1 - 2\,a_2 \right) ) + 
     \left( 1920\,  \right) \,\cosh (2\,\left( a_1 + 2\,a_2 \right) ) - \\
     \left( 256\,  \right) \,\sinh (2\,a_1) - 
     \left( 384\,  \right) \,\sinh (2\,\left( a_1 - 2\,a_2 \right) ) - \\
     \left( 1152\,  \right) \,\sinh (2\,\left( a_1 + 2\,a_2 \right) ) \big)\;;
\end{eqnarray*}
and 
\begin{eqnarray*}
{\bf c}\;:=\;
202 - 232\,\cosh (4\,a_1) + \cosh (8\,a_1) + 
  2\,\cosh (4\,\left( a_1 - a_2 \right) ) - 232\,\cosh (4\,a_2) + \\
  \cosh (8\,a_2) + 2\,\cosh (4\,\left( a_1 + a_2 \right) ) + 
  64\,\sinh (4\,a_1) + 64\,\sinh (4\,a_2)\;.
\end{eqnarray*}
Let $C$ be a positive constant and set ${\bf m}_C:=e^{-i S}\,\left(C+{\bf \Delta}^2\right)\,e^{i S}$. On then has $|{\bf m}_C|\geq|\CQ_4+\CQ_2+{\bf c} +C|$.
Let ${\Xi}={\Xi}(a_1,a_2)$ be a continuous function. The integral $\int_{\S\times\S}\frac{\Xi}{{\bf m}_C}$ exists provided 
the following integral
\begin{equation*}
{\bf {\mathcal I}}\;:=\;\int_{\S\times\S}\frac{{\Xi}}{\CQ_4+\CQ_2+{\bf c} +C}
\end{equation*}
does.

\noindent In order to establish the existence of the above integral ${\mathcal I}$, we first
 observe that the quadratic form
$A:=\left(
\begin{array}{cc}
 \cosh[4 a_2] & \sinh[2 (a_1 + a_2)]\\
 \sinh[2 (a_1 + a_2)]& \cosh[4 a_1] 
\end{array}
\right)$ is positive definite for all $(a_1,a_2)=:a$. Indeed  its determinant equals $[\cosh(2(a_2-a_1))]^2$.

\noindent We then analyse the behaviour of the quadratic form $B$ defined as $\CQ_2=:{}^\tau\ell B\ell$
with $\ell:=(\ell_1,\ell_2)$. One observes that its trace is strictly positive for large values of $|a|$ as
its dominant term is greater than $\cosh(8a_1)+\cosh(a_2)$. Also for large $|a|$, its determinant behaves as $128(\cosh(12a_1-4a_2)+\cosh(4(a_1-3a_2))+2\cosh(8(a_1-a_2)))= 128\times 4(\cosh(8(a_1-a_2))\cosh^2(2(a_1+a_2)))$. Therefore $B$ is positive definite for large values of $|a|$, say for $|a|>{\bf a}_0>0$.

\noindent Note that for large values of $|a|$, say $|a|>{\bf a }_0$ as well, the independent term ${\bf c}$ behaves (at least)
as $\cosh(8a_1)+\cosh(8a_2)=2\cosh(4(a_1+a_2))\cosh(4(a_1-a_2))$.

\noindent Now, one has 
\begin{eqnarray*}
{\mathcal I}&=&\int \frac{\Xi}{({}^\tau\ell A \ell)^2+{}^\tau\ell B \ell+{\bf c} +C}\;,
\end{eqnarray*}
which equals:
\begin{equation*}
\int \frac{\Xi}{({}^\tau\ell \Lambda \ell)^2+{}^\tau\ell U B {}^\tau U \ell+{\bf c} +C}\;,
\end{equation*}
after changing the variables following $\ell\to U\ell$ where $U\in SO(2)$ is such that
${}^\tau U\Lambda U=A$ with $\Lambda$ diagonal (and positive definite).
Hence setting $x:=\Lambda^{\frac{1}{2}}\ell$, one gets:
\begin{equation*}
{\mathcal I}=\int\frac{\Xi}{\sqrt{\det(A)}\,(x^4+{}^\tau x\Lambda^{-\frac{1}{2}}U B {}^\tau U\Lambda^{-\frac{1}{2}}x+{\bf c}+C)}\,{\rm d}x\,{\rm d}a\;.
\end{equation*}
Passing in polar coordinates $x=:r\,e^{i\theta}$, and setting 
\begin{equation*}
\beta\;:=\;{}^\tau e^{i\theta}\Lambda^{-\frac{1}{2}}U\, B\, {}^\tau U\Lambda^{-\frac{1}{2}}e^{i\theta}\;,
\end{equation*}
one obtains
\begin{equation*}
{\mathcal I}=\int_{\R^2}\frac{\Xi}{\sqrt{\det(A)}}\int_0^{2\pi}\int_0^\infty\frac{r}{r^4+\beta\,r^2+{\bf c}+C}\,{\rm d}r\,{\rm d}\theta\,{\rm d}a\;,
\end{equation*}
or, after $r\to r^2$,
\begin{equation*}
{\mathcal I}=\frac{1}{2}\,\int_{\R^2}\frac{\Xi}{\sqrt{\det(A)}}\int_0^{2\pi}\int_0^\infty\frac{1}{r^2+\beta\,r+{\bf c}+C}\,{\rm d}r\,{\rm d}\theta\,{\rm d}a\;.
\end{equation*} 
Now, for $|a|\geq{\bf a}_0$, $\beta$ is strictly positive. Hence $r^2+\beta r+{\bf c}+C>r^2+{\bf c}\sim r^2+2\cosh(4(a_1+a_2))\cosh(4(a_1-a_2))$. Therefore the quantity
\begin{equation*}
{\mathcal I}_+\;:=\;\frac{1}{2}\,\int_{|a|>{\bf a}_0}\frac{\Xi}{\sqrt{\det(A)}}\int_0^{2\pi}\int_0^\infty\frac{1}{r^2+\beta\,r+{\bf c}+C}\,{\rm d}r\,{\rm d}\theta\,{\rm d}a\;,
\end{equation*}
is bounded by
\begin{equation*}
\gamma\,\int_{|a|>{\bf a}_0}\frac{\Xi}{\sqrt{2\det(A)\cosh(4(a_1+a_2))\cosh(4(a_1-a_2))}}\,{\rm d}a\,\int_0^\infty\frac{1}{s^2+1}\,{\rm d}s\;,
\end{equation*}
where $\gamma$ is a positive constant.

\noindent Regarding the quantity:
\begin{equation*}
{\mathcal I}_-\;:=\;\frac{1}{2}\,\int_{|a|\leq{\bf a}_0}\frac{\Xi}{\sqrt{\det(A)}}\int_0^{2\pi}\int_0^\infty\frac{1}{r^2+\beta\,r+{\bf c}+C}\,{\rm d}r\,{\rm d}\theta\,{\rm d}a\;,
\end{equation*}
one may adapt the constant $C$ in such a way that $\beta r + {\bf c} +C\geq1$ for all
$(r,\theta,a)$ in the compact set $[0,1]\times[0,2\pi]\times\{|a|\leq{\bf a}_0\}$. Which yields 
a finite quantity ${\mathcal I}_-$ for any data of a continuous function $\Xi$.
\EPf

\subsection{Deformed products}
 We denote by $\CD_E$, $\CS_E$ and
$\CS'_E$ respectively the space of $E$-valued smooth compactly supported functions
on  $\S$, the space of $E$-valued smooth rapidly decreasing functions
on $\S$ and the space of $E$-valued tempered distributions on $\S$. By nuclearity and denoting by
$\hat{\otimes}_\pi$ the projective tensor product, one has 
the following isomorphisms \cite{Tr}: $\CD_E\simeq\CD\hat{\otimes}_\pi E\;,\;\CS_E\simeq\CS\hat{\otimes}_\pi E$
and $\CS'_E\simeq\CS'\hat{\otimes}_\pi E$. In particular, the partial Fourier transform $\CF$ induces
the following automorphisms: $\CF\hat{\otimes}_\pi\id_E:\CS_E\to\CS_E$ and $\CF\hat{\otimes}_\pi\id_E:\CS'_E\to\CS'_E$.
\begin{lem}
The linear map:
\begin{equation*}
\CS(\R)\longrightarrow\CS'(\R):\varphi\mapsto\varphi\circ\arcsinh
\end{equation*}
is continuous w.r.t. the strong topology on $\CS'(\R)$.
\end{lem}
\Pf
First observe that the linear map
\begin{equation*}
\CS\to\CS:\varphi\mapsto[\tilde{\varphi}:t\mapsto\cosh(t)\varphi(\sinh(t))]
\end{equation*}
is bounded.
Now, assume $\{\phi_n\}$ is a sequence of Schwartz functions that tends to zero in the $\CS$-topology and 
consider the corresponding sequence $\{T_n:=\phi_n\circ\arcsinh\}$ in $\CS'(\R)$. One then has
$<T_n,\varphi>=\int\phi_n\tilde{\varphi}=<\phi_n,\tilde{\varphi}>$, which tends to zero
uniformly for $\varphi$ running in any bounded subset of $\CS(\R)$ as the natural inclusion
$\CS(\R)$ in its strong dual $\CS'(\R)$ is continuous.
\EPf
By nuclearity and for every Schwartz multiplier $\Theta\in\CO_M(\R)$ (\cite{Tr} p. 275), one gets a continuous linear injection:
\begin{equation*}
U^\Theta\hat{\otimes}_\pi\id_E:\CS_E\longrightarrow\CS'_E\;.
\end{equation*}
 We then set
\begin{equation*}
\CE^\Theta_E\;:=\;U^\Theta\hat{\otimes}_\pi\id_E(\,\CS_E\,)\;.
\end{equation*}
At last, we denote by $\Box$ the differential operator on $C^\infty(\S\times\S)$ defined as $\Box F:=(C+{\bf \Delta}^2)\left(\frac{1}{{\bf m}_C}\,F\right)$.  We observe
\begin{lem}\label{PSD}
\begin{enumerate}
\item[(i)] Let $\Theta\in C^\infty(\R)$ be nowhere vanishing and satisfying the property that there exists $N\in\R$ such that one may find $C>0$ with $|\frac{{\rm d}^r}{{\rm d}t^r}\Theta(t)|\leq C(1+|t|)^{N-r}$. Then the function (see Theorem \ref{CONTRACTED})
\begin{equation*}
{\bf \Theta}(x_1,x_2)\;:=\;\frac{{\bf \Xi}(a_1)\,{\bf \Xi}(-a_2)}{{\bf \Xi}(a_1-a_2)}
\end{equation*}
belongs to $\CB^1_\C(\S\times\S)$.
\item[(ii)] The function
\begin{equation*}
{\bf A}_{\mbox{\rm can}}(x_1,x_2)\;:=\;\sqrt{\left|{\mbox{\rm Jac}}_{\Phi}(e,x_1,x_2)\right|}
\end{equation*}
is an amplitude on $\S\times\S$ adapted to $S_{\mbox{\rm can}}$.
\end{enumerate}
\end{lem}
We now consider a strongly continuous isometric action $\alpha$ of $\S$ on a Fr\'echet algebra $(\A.\mu_\A)$ and denote
by $\A_\infty$ the associated space of smooth vectors. 
\begin{cor} Let $\Theta\in C^\infty(\R)$ be as in Lemma \ref{PSD}.
For all $a,b\in\A_\infty$, the following oscillatory integral 
\begin{equation}\label{UDF}
a\star^\Theta_\theta b:=\frac{1}{\theta^2}\tilde{\int}\,e^{\frac{i}{\theta}S_{\mbox{\rm can}}}\,{\bf \Theta}\,{\bf A}_{\mbox{\rm can}}\,\mu_\A(\alpha(a)\otimes\alpha(b))
\end{equation}
is well defined as an element of $\A_\infty$.
\end{cor}
\begin{thm}
The pair $(\A_\infty,\star^\Theta_\theta)$ is an associative topological algebra.
\end{thm}
\Pf
The space $\CE_\theta(\S)\hat{\otimes}\A_\infty$ naturally inherits
a structure of associative algebra denoted by $\hat{\star}$ from the one on $\CE_\theta(\S)$. Consider $a,b,c\in\A_\infty$, then we have
\begin{eqnarray*}
(a\star b)\star c&=&\tilde{\int}K(e,g_1,g_2)\left(\alpha_{g_1}\tilde{\int}K(e,h_1,h_2)\alpha_{h_1}(a)\,\alpha_{h_2}(b)\;{\rm d}h_1\,{\rm d}h_2\right)\,\alpha_{g_2}(c)\;{\rm d}g_1\,{\rm d}g_2\\
&=&{\int} e^{iS_e(g_1,g_2)}\Box_{(g_1,g_2)}\left(\alpha_{g_1}\left(\tilde{\int}K(e,h_1,h_2)\alpha_{h_1}(a)\,\alpha_{h_2}(b)\;{\rm d}h_1\,{\rm d}h_2\right)\,\alpha_{g_2}(c)\right)
\;{\rm d}g_1\,{\rm d}g_2\\
&=&{\int} e^{iS_e(g_1,g_2)}\Box_{(g_1,g_2)}\left(\alpha_{g_1}\left({\int}\,e^{iS_e(h_1,h_2)}\Box_{(h_1,h_2)}\left(\alpha_{h_1}(a)\,\alpha_{h_2}(b)\right)\;{\rm d}h_1\,{\rm d}h_2\right)\,\alpha_{g_2}(c)\right)
\;{\rm d}g_1\,{\rm d}g_2\\
&=&\int\int e^{i(S_e(g_1,g_2)+S_e(h_1,h_2))}\Box_{(g_1,g_2)}\Box_{(h_1,h_2)}\left(\alpha_{g_1}\left(\alpha_{h_1}(a)\,\alpha_{h_2}(b)\right)\right)\,\alpha_{g_2}(c)
\;{\rm d}h_1\,{\rm d}h_2\;{\rm d}g_1\,{\rm d}g_2\\
&=&\int\int e^{i(S_e(g_1,g_2)+S_e(h_1,h_2))}\Box_{(h_1,h_2)}\Box_{(g_1,g_2)}\left(\alpha_{g_1h_1}(a)\,\alpha_{g_1h_2}(b)\right)\,\alpha_{g_2}(c)
\;{\rm d}h_1\,{\rm d}h_2\;{\rm d}g_1\,{\rm d}g_2\\
&=&\int\int e^{i(S_e(g_1,g_2)+S(g_1,h_1,h_2))}\Box_{(g_1^{-1}h_1,g_1^{-1}h_2)}\Box_{(g_1,g_2)}\left(\alpha_{h_1}(a)\,\alpha_{h_2}(b)\,\alpha_{g_2}(c)\right)
\;{\rm d}h_1\,{\rm d}h_2\;{\rm d}g_1\,{\rm d}g_2\\
&=&\lim_n\left[\int\int e^{i(S_e(g_1,g_2)+S(g_1,h_1,h_2))}\Box_{(g_1^{-1}h_1,g_1^{-1}h_2)}\Box_{(g_1,g_2)}
\left([\alpha(a)]_n(h_1)\,[\alpha(b)]_n(h_2)\,[\alpha(c)]_n(g_2)\right)
\right]\\
&=&\lim_n\left[\int\int e^{i(S_e(g_1,g_2)+S_e(h_1,h_2))}\Box_{(h_1,h_2)}\Box_{(g_1,g_2)}
\left([\alpha(a)]_n(g_1h_1)\,[\alpha(b)]_n(g_1h_2)\,[\alpha(c)]_n(g_2)\right)
\right]\\
&=&\lim_n\left[\int\int e^{i(S_e(g_1,g_2)+S_e(h_1,h_2))}\Box_{(g_1,g_2)}\Box_{(h_1,h_2)}
\left([\alpha(a)]_n(g_1h_1)\,[\alpha(b)]_n(g_1h_2)\,[\alpha(c)]_n(g_2)\right)
\right]\\
&=&\lim_n\left[\int e^{iS_e(g_1,g_2)}\Box_{(g_1,g_2)}\int \left(e^{iS_e(h_1,h_2)}\Box_{(h_1,h_2)}
\left([\alpha(a)]_n(g_1h_1)\,[\alpha(b)]_n(g_1h_2)\,\right)\right)\,[\alpha(c)]_n(g_2)
\right]\\
&=&\lim_n\left[\int e^{iS_e(g_1,g_2)}\Box_{(g_1,g_2)}\int \left(e^{iS_e(g_1,h_1,h_2)}\Box_{(h_1,h_2)}
\left([\alpha(a)]_n(h_1)\,[\alpha(b)]_n(h_2)\,\right)\right)\,[\alpha(c)]_n(g_2)
\right]\\
&=&\lim_n\left[\int e^{iS_e(g_1,g_2)}\Box_{(g_1,g_2)}
\left([\alpha(a)]_n\hat{\star}\,[\alpha(b)]_n\,\right)(g_1)\,[\alpha(c)]_n(g_2)
\right]\\
&=&\lim_n\left[\tilde{\int} \,K(e,g_1,g_2)\,
\left([\alpha(a)]_n\hat{\star}\,[\alpha(b)]_n\,\right)(g_1)\,[\alpha(c)]_n(g_2)
\right]\\
&=&\lim_n\left(
\left([\alpha(a)]_n\hat{\star}\,[\alpha(b)]_n\,\right)\hat{\star}[\alpha(c)]_n
\right)(e)\\
&=&\lim_n\left(
[\alpha(a)]_n\hat{\star}\,\left([\alpha(b)]_n\,\hat{\star}[\alpha(c)]_n\right)
\right)(e)\\
&=& a\star(b\star c)\;.
\end{eqnarray*}
\EPf

\end{document}